\documentclass[11pt]{article}
\usepackage[a4paper,margin=25mm]{geometry}
\usepackage{amsmath,amssymb,amsthm,booktabs,tikz,microtype,graphicx}
\usetikzlibrary{arrows.meta}
\usepackage[colorlinks=true,linkcolor=blue!45!black,citecolor=blue!45!black,urlcolor=blue!45!black]{hyperref}
\usepackage{authblk}
\newtheorem{theorem}{Theorem}[section]
\newtheorem{lemma}[theorem]{Lemma}
\newtheorem{proposition}[theorem]{Proposition}
\newtheorem{corollary}[theorem]{Corollary}
\theoremstyle{definition}\newtheorem{definition}[theorem]{Definition}
\theoremstyle{remark}\newtheorem{remark}[theorem]{Remark}

\newcommand{\din}{d^-}
\newcommand{\topsub}{\preccurlyeq_{\mathrm{top}}}
\newtheorem{example}[theorem]{Example}

\title{Longitudinal Subdivision Invariance of Periodic Directed Strips and Its Application to Elongated Octagonal Strips}
\author{Dan Zhu\thanks{E-mail: zd2609477480@163.com}}

\affil{School of Mathematical Sciences, Shenzhen University,
Shenzhen 518060, China}

\date{}

\begin{document}

\maketitle

\vspace{4pt}

\begin{abstract}
We prove that longitudinal refinements of finite-width periodic directed strips preserve subdivision containment. The refinement lengths may be arbitrary positive integers, without a uniform upper bound, and no matching condition is imposed on the horizontal arcs. We construct subdivision models in both directions that preserve each designated rail as a whole. The equivalence also holds when the designated rail images are required to lie in a prescribed end. A nonperiodic example  shows that reverse containment may fail without periodicity. As an application, we obtain an elongated octagonal strip $O_n$ by subdividing selected longitudinal arcs of the directed hexagonal strip $H_n$. We show that every bounded face of $O_n$ is an octagon for $n\geq 2$ and that $O_n$ and $H_n$ are subdivision-equivalent. Consequently, the known exact end in-degree thresholds and forced-width formula for $H_n$ also hold for $O_n$.
\end{abstract}

\section{Introduction}
Halin's grid theorem \cite{Halin} establishes a fundamental connection between the structure of ends and the existence of grid-like subgraphs in infinite graphs. Hamann and Heuer \cite{HH} established directed analogues of Halin's grid theorem and showed that sufficiently large end in-degree forces subdivisions of directed hexagonal grids of prescribed finite width. This raises the natural question of determining the smallest end in-degree that guarantees such a subdivision. Zhu and Lyu \cite{ZL} subsequently determined the exact threshold for each width, including the exceptional cases of widths one, two, and three. They also obtained the corresponding forced-width formula. These results motivate the study of related directed strips.

We consider longitudinal refinements of finite-width periodic directed strips. Such a refinement replaces each longitudinal arc by a finite directed path while leaving the horizontal arcs unchanged. Every longitudinal refinement contains a subdivision of the original strip, but the reverse containment may fail without periodicity.
 We investigate whether periodicity of the horizontal arc pattern guarantees subdivision equivalence.

Our main result shows that every longitudinal refinement of a finite-width periodic directed strip is subdivision-equivalent to the original strip. The refinement lengths need not be uniformly bounded, and horizontal arcs may share endpoints. The subdivision models preserve each designated rail as a whole, so the equivalence also holds when the rail images are required to lie in a specified end. We also give a nonperiodic example for which reverse containment fails.

We apply this result to an elongated octagonal strip obtained by subdividing selected longitudinal arcs of the directed hexagonal strip. We verify its octagonal face structure and prove subdivision equivalence with the hexagonal strip. The exact end in-degree thresholds and forced-width formula for $O_n$ then follow from the corresponding results for $H_n$ in \cite{ZL}. We also discuss the limits of the refinement theorem and the distinction between these elongated strips and standard square-octagon lattices.

\section{Preliminaries}\label{sec:prelim}
All digraphs considered here have no loops or parallel arcs. Oppositely directed arcs between two vertices are allowed. Let $\mathbb{N}_0=\{0,1,2,\ldots\}$. All directed paths are finite unless stated otherwise. A ray is a one-way infinite directed path. An anti-ray is a one-way infinite path whose arcs are directed towards its terminal vertex. When defining the relation between rays and anti-rays, we allow a path of length zero.

For any two rays or anti-rays $R,S$ in a digraph $D$, a directed $R$--$S$ path is a directed path with initial vertex in $V(R)$, terminal vertex in $V(S)$, and no internal vertex in $V(R)\cup V(S)$. Write $R\le S$ if there are infinitely many pairwise vertex-disjoint directed $R$--$S$ paths, and write $R\sim S$ if $R\le S$ and $S\le R$. Ends are the equivalence classes of rays and anti-rays under the corresponding relation \cite{Zuther,HH}. We consider only ends containing rays. Following Hamann and Heuer \cite{HHdegree}, the in-degree $\din(\omega)$ of an end $\omega$ is the maximum number of pairwise vertex-disjoint rays in $\omega$, with $d^-(\omega)=\infty$ if arbitrarily large finite families of such rays exist.

A directed subdivision model of a digraph $F$ in a digraph $D$ assigns distinct vertices of $D$ to the vertices of $F$ and replaces each arc of $F$ by a finite directed path of positive length between the corresponding assigned vertices. The replacement paths have pairwise disjoint interiors, and no internal vertex is an assigned vertex. We write $F\preceq_{\mathrm{top}}D$ if such a model exists. Two digraphs $F$ and $G$ are called subdivision-equivalent if $F\preceq_{\mathrm{top}}G$ and $G\preceq_{\mathrm{top}}F$. The model need not be induced, so unused arcs of $D$ are irrelevant.

If $F$ has specified pairwise vertex-disjoint rays, we call their images in a subdivision model the designated rail images. Each image is a directed ray obtained by concatenating the replacement paths of the corresponding rail. This follows from the distinctness of the assigned vertices and the disjointness of the replacement interiors. The terminology also applies to width one and does not require every rail vertex to have degree three.

\begin{lemma}\label{lem:composition}
If $F\topsub G$ and $G\topsub D$, then $F\topsub D$. Moreover, if the model of $F$ in $G$ maps each designated rail of $F$ onto an entire designated rail of $G$, then the composed model maps that rail onto the corresponding designated rail image in $D$.
\end{lemma}
\begin{proof}
For each arc of $F$, consider its replacement path in $G$. Replace each arc of this path by its replacement path in  $D$ and concatenate the resulting paths. The concatenation is finite and directed. Since the original path in $G$ has no repeated vertices, and the replacement paths in $D$ have pairwise disjoint interiors and avoid all assigned
vertices of $G$, the concatenation is a path. The same properties ensure that the composed paths have pairwise
disjoint interiors and that no internal vertex is an assigned vertex of $F$. Thus they define a subdivision
model of $F$ in $D$.

For the second assertion, each designated rail of $F$ is represented by the entire corresponding rail of $G$. The composition therefore uses every replacement path of that rail in $D$, in the same order. Hence its image is exactly the corresponding designated rail image in $D$.
\end{proof}

\section{Longitudinal refinements of periodic strips}\label{sec:periodic}
Fix integers $n,p\ge1$, and write $[n]=\{1,\ldots,n\}$. For each $r\in\{0,\ldots,p-1\}$, let
\[
A_r\subseteq
\{(a,b)\in[n]\times[n]:a\neq b\}.
\]
The set $A_r$ specifies the horizontal arcs at levels
$j\equiv r\pmod p$.
Write
\[
\mathcal A=(A_0,A_1,\ldots,A_{p-1}).
\]
Some of the sets $A_r$ may be empty. No matching condition is imposed. Distinct arcs in $A_r$ may share endpoints, and both $(a,b)$ and $(b,a)$ may belong to $A_r$.

\begin{definition}\label{def:periodic}
The periodic directed strip $P=P(n,\mathcal A)$ has vertices $v_j^i$ for $i\in[n]$ and $j\in\mathbb N_0$, longitudinal arcs
\[
v_j^iv_{j+1}^i,
\] and horizontal arcs 
\[
v_j^av_j^b
\qquad
\text{for every }(a,b)\in A_{j\bmod p}.
\]
Its designated rails are 
\[
v_0^iv_1^iv_2^i\cdots,
\qquad i\in[n].
\]

Let $\ell=(\ell_{i,j})_{i\in[n],j\in\mathbb N_0}$ be any array of positive integers. The longitudinal refinement $P[\ell]$ replaces the longitudinal arc $v_j^iv_{j+1}^i$ by a directed path of length $\ell_{i,j}$, retaining all horizontal arcs. Every inserted vertex is new, and distinct replacement paths have disjoint interiors. Its designated rails are the refined longitudinal rays.
\end{definition}

\begin{theorem}\label{thm:refinement}
For every $n,p\geq 1$, every tuple $\mathcal A=(A_0,\ldots,A_{p-1})$ of horizontal arc sets as in Definition \ref{def:periodic}, and every positive integer array $\ell$,
\[
P(n,\mathcal A)
\preceq_{\mathrm{top}}
P(n,\mathcal A)[\ell]
\]
and
\[
P(n,\mathcal A)[\ell]
\preceq_{\mathrm{top}}
P(n,\mathcal A).
\]
Both models represent each designated rail by the entire corresponding designated rail. Consequently, for every digraph $D$,
\[
P(n,\mathcal A)\preceq_{\mathrm{top}}D
\quad\Longleftrightarrow\quad
P(n,\mathcal A)[\ell]\preceq_{\mathrm{top}}D.
\]
The equivalence remains valid with the additional requirement that all designated rail images belong to any specified end $\omega$ of $D$.
\end{theorem}

\begin{proof}
Write $P=P(n,\mathcal A)$ and $Q=P[\ell]$. The definition of $Q$ gives a subdivision model of $P$ in $Q$, in which each designated rail is mapped onto the entire corresponding rail of $Q$.

For the reverse direction, set
\[
L_j=\max_{1\le i\le n}\ell_{i,j}.
\]
Since $n$ is finite, $L_j$ is finite for each $j$. Set $t_0=0$ and define
\begin{equation}\label{eq:t}
t_{j+1}=t_j+1+p\left\lceil\frac{L_j-1}{p}\right\rceil.
\end{equation}
Then $t_j\equiv j\pmod p$, and
\begin{equation}\label{eq:gap}
L_j\le t_{j+1}-t_j\le L_j+p-1.
\end{equation}
In particular, the sequence $(t_j)$ is strictly increasing and unbounded.

Map each retained vertex $v_j^i$ of $Q$ to $v_{t_j}^i$ in $P$. Denote the inserted vertices on the refined $v_j^i$--$v_{j+1}^i$ path, in order, by
\[
u_{i,j,1},\ldots,u_{i,j,\ell_{i,j}-1}.
\]
For $1\le k<\ell_{i,j}$, map $u_{i,j,k}$ to $v_{t_j+k}^i$. By \eqref{eq:gap}, these vertices lie strictly between $v_{t_j}^i$ and $v_{t_{j+1}}^i$ on rail $i$. Indeed, $k<\ell_{i,j}\le L_j\le t_{j+1}-t_j$, so $t_j<t_j+k<t_{j+1}$.

Represent each arc of the refined path by the longitudinal segment between its assigned endpoints. If $\ell_{i,j}>1$, all but its last arc are represented by single longitudinal arcs, and its last arc by the remaining segment to $v_{t_{j+1}}^i$. If $\ell_{i,j}=1$, its only arc is represented by the whole segment. All these segments have positive finite length.

For a horizontal arc $v_j^av_j^b$ of $Q$, where $(a,b)\in A_{j\bmod p}$, use the single arc $v_{t_j}^av_{t_j}^b$ of $P$. It exists because $t_j\equiv j\pmod p$, and its direction is unchanged. Different horizontal arcs may share assigned endpoints, but this does not violate the subdivision conditions. Each such arc is represented by a single arc and therefore has no internal vertices. Thus, any common endpoint of these replacement paths is the image of a vertex of $Q$, and their interiors remain disjoint.

The vertex assignment is injective. Assigned vertices on different rails are distinct, and on one rail the intervals between successive $t_j$ are disjoint apart from their shared endpoints. Within each interval, the assigned vertices occur in strictly increasing order. Thus the longitudinal replacement paths have pairwise disjoint interiors and avoid every other assigned vertex. Horizontal replacement paths have no internal vertices and meet longitudinal paths only at the prescribed endpoints. Discarding all unused horizontal arcs therefore gives a subdivision
model of $Q$ in $P$.

Since $t_0=0$ and $t_j\to\infty$, the longitudinal segments used on rail $i$ cover its entire ray. Both models preserve the whole designated rails. The containment equivalence follows from Lemma~\ref{lem:composition}, and its final assertion follows because composition retains exactly the same rail images in $D$.
\end{proof}

\begin{remark}\label{rem:refinement}
Theorem~\ref{thm:refinement} allows arbitrary finite longitudinal
subdivision lengths without a uniform upper bound.
No matching condition is imposed on the horizontal
arc sets. The theorem does not concern arbitrary
subdivisions of horizontal arcs.
\end{remark}

\begin{example}\label{ex:nonperiodic}
Let $G$ have vertex set
\[
V(G)=\{x_j^i:i\in\{1,2\},\ j\in\mathbb N_0\},
\]
with longitudinal arcs
\[
x_j^i\to x_{j+1}^i
\qquad (i\in\{1,2\},\ j\in\mathbb N_0).
\]
Add the two horizontal arcs
\[
x_1^1\to x_1^2,
\qquad
x_2^2\to x_2^1.
\]
Since horizontal arcs occur only at levels $1$ and $2$, the horizontal arc pattern is not periodic.

Let $G'$ be obtained from $G$ by subdividing each arc $x_1^i\to x_2^i$, $i\in\{1,2\}$, once.
The underlying undirected graph of $G$ has exactly one cycle of length four, whereas that of $G'$ has exactly one cycle of length six. If $G'\preceq_{\mathrm{top}}G$, then the $6$-cycle in $G'$ would give a subdivision of a $6$-cycle in the underlying undirected graph of $G$. Such a subdivision contains a cycle of length at least six, whereas the underlying undirected graph of $G$ has no cycle other than a $4$-cycle. Hence
\[
G'\not\preceq_{\mathrm{top}}G.
\]
Thus the reverse containment can fail when the horizontal arc pattern is nonperiodic.
\end{example}

\section{Elongated octagonal strips and end in-degree thresholds}\label{sec:octagonal}
For $n\ge1$, let $H_n$ be the digraph with vertex set
\[
V(H_n)=\{x_j^i\mid i\in[n],\ j\in\mathbb{N}_0\}
\]
and longitudinal arcs $x_j^ix_{j+1}^i$ for all $i\in[n]$ and $j\in\mathbb{N}_0$. 

For $r\in\{0,1,2,3\}$, define 
\[
\begin{aligned}
M_0
&=\{(i+1,i):1\leq i<n,\ i\text{ even}\},\\
M_1
&=\{(i,i+1):1\leq i<n,\ i\text{ odd}\},\\
M_2
&=\{(i,i+1):1\leq i<n,\ i\text{ even}\},\\
M_3
&=\{(i+1,i):1\leq i<n,\ i\text{ odd}\}.
\end{aligned}
\]
At each level $j\in\mathbb N_0$, the horizontal arcs of $H_n$ are precisely
\[
x_j^a\to x_j^b
\qquad\text{for all }(a,b)\in M_{j\bmod 4}.
\]
There are no other horizontal arcs. Each $M_r$ is a directed matching. Hence, in the notation of Definition~\ref{def:periodic},
\[
H_n=P(n,\mathcal M),
\qquad
\mathcal M=(M_0,M_1,M_2,M_3).
\]

\begin{definition}\label{def:O}
For $n\ge 1$, let $O_n=H_n[\ell]$, where
\[
\ell_{i,j}=\begin{cases}1,&j\text{ even},\\2,&j\text{ odd}.\end{cases}
\]
Equivalently, $O_n$ is obtained from $H_n$ by subdividing each longitudinal arc $x_{2t+1}^i x_{2t+2}^i$ exactly once
for $i\in[n]$ and $t\in\mathbb{N}_0$, while horizontal arcs and their directions remain unchanged. For $n=1$, $O_1$ is a directed ray.
\end{definition}

We use the natural plane embedding in which the longitudinal rays are drawn as parallel vertical rays and the horizontal edges join vertices at the same level.

\begin{proposition}\label{prop:faces}
For $n\ge2$, every bounded face of the underlying undirected graph of $O_n$ is bounded by an $8$-cycle in this  plane embedding.
\end{proposition}
\begin{proof}
By the definitions of $M_0,M_1,M_2,M_3$, between two consecutive rails of $H_n$, successive horizontal edges occur at levels $j$ and $j+2$ of the same parity. On each rail, the segment between these levels consists of two edges, exactly one of which starts at an odd level. By Definition~\ref{def:O}, this edge is subdivided once in $O_n$, so each longitudinal side contains three edges. Together with the two horizontal edges, they form a cycle of length $3+1+3+1=8$. Since within each strip between adjacent rails the horizontal edges separate exactly the regions between consecutive horizontal edges, these are all bounded faces. Longitudinal subdivision does not create additional faces.
\end{proof}

\begin{figure}[ht]
\centering
\resizebox{0.85\textwidth}{!}{%
\begin{tikzpicture}[>=Stealth,scale=1,vertex/.style={circle,fill=black,inner sep=1.15pt}]
\begin{scope}[x=.8cm,y=.46cm]
\foreach \i in {0,...,3}{
 \foreach \j in {0,...,7}{\pgfmathtruncatemacro{\k}{\j+1}\draw[->,thin] (\i,-\j)--(\i,-\k);}
 \foreach \j in {0,...,8}{\node[vertex] at (\i,-\j) {};}
 \draw[densely dotted] (\i,-8)--(\i,-8.6);
}
\foreach \i in {0,...,2}{\foreach \j in {0,...,8}{
 \pgfmathtruncatemacro{\r}{mod(\j,4)}
 \pgfmathtruncatemacro{\railparity}{mod(\i,2)}
 \pgfmathtruncatemacro{\next}{\i+1}
 \ifnum\railparity=0
  \ifnum\r=1 \draw[->] (\i,-\j)--(\next,-\j);\fi
  \ifnum\r=3 \draw[<-] (\i,-\j)--(\next,-\j);\fi
 \else
  \ifnum\r=2 \draw[->] (\i,-\j)--(\next,-\j);\fi
  \ifnum\r=0 \draw[<-] (\i,-\j)--(\next,-\j);\fi
 \fi
}}
\node at (1.5,1) {$H_4$};
\end{scope}
\begin{scope}[shift={(6,0)},x=.8cm,y=.307cm]
\foreach \i in {0,...,3}{
 \foreach \j in {0,...,11}{\pgfmathtruncatemacro{\k}{\j+1}\draw[->,thin] (\i,-\j)--(\i,-\k);}
 \foreach \j in {0,...,12}{
  \pgfmathtruncatemacro{\r}{mod(\j,3)}
  \ifnum\r=2 \draw[red!70!black,fill=white] (\i,-\j) circle[radius=1.7pt];
  \else \node[vertex] at (\i,-\j) {};\fi
 }
 \draw[densely dotted] (\i,-12)--(\i,-12.9);
}
\foreach \i in {0,...,2}{\foreach \j in {0,...,12}{
 \pgfmathtruncatemacro{\r}{mod(\j,6)}
 \pgfmathtruncatemacro{\railparity}{mod(\i,2)}
 \pgfmathtruncatemacro{\next}{\i+1}
 \ifnum\railparity=0
  \ifnum\r=1 \draw[->] (\i,-\j)--(\next,-\j);\fi
  \ifnum\r=4 \draw[<-] (\i,-\j)--(\next,-\j);\fi
 \else
  \ifnum\r=3 \draw[->] (\i,-\j)--(\next,-\j);\fi
  \ifnum\r=0 \draw[<-] (\i,-\j)--(\next,-\j);\fi
 \fi
}}
\node at (1.5,1.5) {$O_4$};
\end{scope}
\end{tikzpicture}}
\caption{Initial portions of $H_4$ and $O_4$, with rays directed downwards. Open red vertices are inserted on longitudinal arcs. The two drawings have different longitudinal scales.}
\label{fig:strips}
\end{figure}
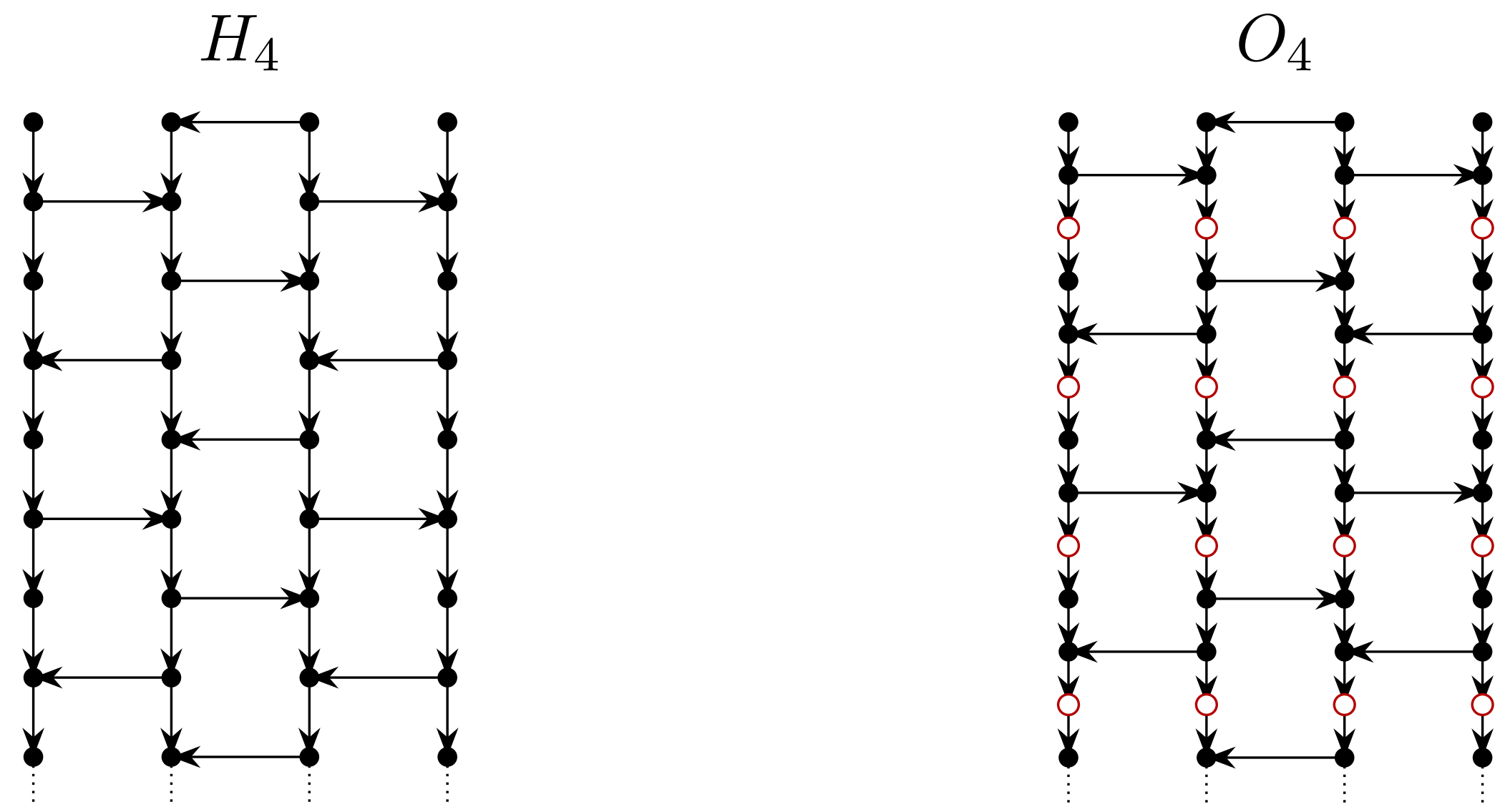

Figure~\ref{fig:strips} illustrates initial portions of $H_4$ and $O_4$.

\begin{corollary}\label{cor:oct}
For every $n\ge1$ and every digraph $D$, $H_n\topsub D$ if and only if $O_n\topsub D$. Moreover, for every end $\omega$ of $D$, $H_n$ has a subdivision model in $D$ whose designated rail images lie in $\omega$ if and only if $O_n$ has such a model.
\end{corollary}
\begin{proof}
By Definition~\ref{def:O}, $O_n$ is a longitudinal refinement of $H_n$. Theorem~\ref{thm:refinement} gives both the stated equivalence and its end-local version.
\end{proof}

\begin{remark}
For each $i\in[n]$, index the vertices on the $i$-th rail of $O_n$ as $y_0^i,y_1^i,\ldots$. The original vertex $x_k^i$ of $H_n$ appears as $y_{f(k)}^i$, where $f(k)=k+\lfloor k/2\rfloor$. Hence the horizontal arcs of $O_n$ occur only at levels congruent to $0,1,3,4$ modulo $6$.

For $t\in\mathbb N_0$ and $r\in\{0,1,2,3,4,5\}$, set
\[
g(6t+r)=12t+a_r,\qquad
(a_0,a_1,a_2,a_3,a_4,a_5)=(0,1,2,6,7,8).
\]
Map $y_j^i$ to $x_{g(j)}^i$ and represent each longitudinal arc by the corresponding directed longitudinal segment. The map $g$ is strictly increasing, with $g(0)=0$ and $g(j)\to\infty$. 

The horizontal levels congruent to $0,1,3,4$ modulo $6$ are mapped to levels congruent to $0,1,2,3$ modulo $4$,
respectively. By the definitions of $M_0,M_1,M_2,M_3$, the corresponding horizontal arcs have the same directions.

The longitudinal segments have pairwise disjoint interiors and contain no other assigned vertices. Horizontal arcs are represented by single arcs. Hence these paths form a subdivision model of $O_n$ in $H_n$ that preserves every designated rail.
\end{remark}

Let $F=(F_n)_{n\ge 1}$ be a family of finite-width directed strips with designated rails. Define $\kappa_F(n)$ to be the least positive integer $q$ such that every digraph containing an end of in-degree at least $q$ contains a subdivision of $F_n$. Define $\kappa_F^{\mathrm{end}}(n)$ similarly, but require that for every digraph $D$ and every end $\omega$ of $D$ with $d^-(\omega)\ge q$, there is a subdivision model of $F_n$ in $D$ whose designated rail images all belong to $\omega$. If no such positive integer $q$ exists, the corresponding threshold is defined to be $\infty$.

Write $H=(H_n)_{n\geq1}$ and $O=(O_n)_{n\geq1}$. For each $n\geq1$, let
\[
\ell^{(n)}
=(\ell^{(n)}_{i,j})_{i\in[n],\,j\in\mathbb N_0}
\]
be an array of positive integers. Set
\[
\Lambda=(\ell^{(n)})_{n\geq1}
\]
and define
\[
H[\Lambda]
=(H_n[\ell^{(n)}])_{n\geq1}.
\]

\begin{proposition}\label{prop:transfer}
For every $n\ge 1$ and every family of positive integer arrays $\Lambda=(\ell^{(n)})_{n\geq1}$,
\[
\kappa_{H[\Lambda]}(n)=\kappa_H(n),
\qquad
\kappa_{H[\Lambda]}^{\mathrm{end}}(n)
=\kappa_H^{\mathrm{end}}(n).
\]
\end{proposition}
\begin{proof}
By Theorem~\ref{thm:refinement}, $H_n[\ell^{(n)}]\preceq_{\mathrm{top}}D$ if and only if $H_n\preceq_{\mathrm{top}}D$ for every digraph $D$. The equivalence also holds when all designated rail images are required to belong to any specified end of $D$. Hence the defining conditions for each pair of thresholds are satisfied by the same positive integers $q$, so the corresponding thresholds are equal.
\end{proof}

We use the exact threshold results and their end-local forms established by Zhu and Lyu \cite{ZL}.
\begin{theorem}\label{thm:input}(Zhu and Lyu\cite{ZL}).
For every $n\ge 1$,
\[
\kappa_H(n)=\kappa_H^{\mathrm{end}}(n)=K(n),
\qquad
K(n)=\begin{cases}
1,&n=1,\\2,&n=2,\\4,&n=3,\\
\lfloor3n/2\rfloor-1,&n\ge4.
\end{cases}
\]
For each $n\ge2$, there is a digraph $D_n$ with an end $\omega_n$ satisfying $d^-(\omega_n)=K(n)-1$ such that $D_n$ contains no subdivision of $H_n$.
\end{theorem}

\begin{corollary}\label{cor:threshold}
For every $n\ge1$,
\[
\boxed{\kappa_O(n)=\kappa_O^{\mathrm{end}}(n)=K(n).}
\]
In particular, for $n\ge 4$, every digraph $D$ with an end $\omega$ satisfying $d^-(\omega)\ge \lfloor 3n/2\rfloor-1$ contains a subdivision of $O_n$ whose designated rail images all belong to $\omega$. For each $n\ge2$, there exists a digraph with an end of in-degree $K(n)-1$ that contains no subdivision of $O_n$.
\end{corollary}

\begin{proof}
By Corollary~\ref{cor:oct}, the subdivision-containment equivalence between $H_n$ and $O_n$ holds both without restrictions and with all designated rail images required to belong to any specified end. The threshold equalities therefore follow from Theorem~\ref{thm:input}.

For $n\ge 2$, let $D_n$ be the lower-bound example from Theorem~\ref{thm:input}. Since $D_n$ contains no subdivision of $H_n$, Corollary~\ref{cor:oct} implies that it contains no subdivision of $O_n$.
\end{proof}

For a family $F=(F_n)_{n\geq1}$ of finite-width directed strips with designated rails and for every positive integer 
 $d$, let $W_F(d)$ be the largest width $n$, if it exists, such that every digraph $D$ with an end $\omega$ satisfying $d^-(\omega)=d$ contains a subdivision of $F_n$ whose designated rail images all belong to $\omega$.

\begin{corollary}\label{cor:width}
For every positive integer $d$,
\[
W_O(d)=\begin{cases}
1,&d=1,\\
2,&d=2,3,\\
\lfloor2d/3\rfloor+1,&d\ge4.
\end{cases}
\]
Moreover, the same formula holds for every family $F=H[\Lambda]$, where $\Lambda$ is as defined above.
\end{corollary}

\begin{proof}
For every $d,n\geq1$, the end-local subdivision guarantee for $O_n$ is equivalent to that for $H_n$ by Corollary~\ref{cor:oct}. Hence they have the same forced width. The forced-width formula for $H_n$ is established in Section~5 of \cite{ZL}.

Now let $F=H[\Lambda]$ be as stated. For each $n$, Theorem~\ref{thm:refinement} gives the same end-local containment equivalence between $F_n$ and $H_n$. Thus, for every $d\ge 1$, the set of widths guaranteed for $F$ coincides with that for $H$. Consequently,
\[
W_F(d)=W_O(d).
\]
\end{proof}

\section{Scope and concluding remarks}\label{sec:scope}
Theorem~\ref{thm:refinement} relies on the periodicity of the horizontal arc pattern and does not require a plane embedding or restrictions on face lengths. For $O_n$, Proposition~\ref{prop:faces} shows that every bounded face has boundary length eight, but this property is not used to establish subdivision equivalence.

The elongated octagonal strip $O_n$ differs from the standard square--octagon lattice. The former contains vertices of degree two inserted on longitudinal arcs, whereas the latter is cubic and has both square and octagonal faces \cite{KPW}. Thus the present results do not directly apply to general directed square--octagon strips.

A separate question is whether finite end in-degree thresholds exist for more general directed targets.

\begin{proposition}\label{prop:cycle}
If a digraph $F$ contains a directed cycle, then for every positive integer $m$, there exists a digraph $D_m$ with an end of in-degree exactly $m$ such that $D_m$ contains no subdivision of $F$.
\end{proposition}

\begin{proof}
Fix $m\ge 1$ and let $D_m=H_m$. Every arc of $H_m$ either increases the level or lies within a single level. Since no arc decreases the level, any directed cycle would have to lie within one level. This is impossible because the horizontal arcs at each level form a matching. Hence $H_m$ is acyclic.

Let $L_1,\ldots,L_m$ be the designated rails of $H_m$. Consecutive rails are joined in both directions by infinitely many pairwise vertex-disjoint directed paths. Thus all $m$ rails belong to the same end $\omega_m$, giving
\[
d^-(\omega_m)\geq m.
\]

Each level of $H_m$ contains exactly $m$ vertices, and no arc decreases the level. Since each level is finite,
every directed ray has unbounded level. Suppose that $H_m$ contains $m+1$ pairwise vertex-disjoint directed rays.
Choose $k$ larger than the initial levels of all these rays. Each ray must cross the cut between levels $k$ and $k+1$ along a longitudinal arc. Since the rays are pairwise vertex-disjoint, they must use distinct arcs of this cut. However, the cut contains exactly $m$ longitudinal arcs, a contradiction. Thus $H_m$ contains at most $m$ pairwise vertex-disjoint directed rays, and hence
\[
d^-(\omega_m)\leq m.
\]
Together with the reverse inequality, this gives
\[
d^-(\omega_m)=m.
\]

A subdivision of a directed cycle is again a directed cycle. Therefore $D_m=H_m$ has an end of in-degree exactly 
$m$ and contains no subdivision of $F$.
\end{proof}

For example, if the boundary of a square face in a directed square--octagon strip forms a directed cycle, Proposition~\ref{prop:cycle} shows that no finite end in-degree bound guarantees a subdivision of the strip.  Theorem~\ref{thm:refinement} concerns longitudinal refinements of periodic directed strips with unchanged horizontal arcs and does not apply directly to general directed square--octagon strips.

Example~\ref{ex:nonperiodic} shows that reverse subdivision
containment may fail for nonperiodic
horizontal arc patterns.
It remains natural to ask which more general
horizontal arc patterns preserve subdivision
containment under arbitrary longitudinal
refinements. Another question is whether
suitable acyclic orientations of finite-width
square--octagon strips admit analogous
end in-degree threshold results.

\par\vspace{1.5em}
\smallskip
\noindent\textbf{Data availability.}
No data were used for the research described in this paper.

\noindent\textbf{Declaration on the Use of AI.}
AI tools were used to assist with checking proof arguments and improving the manuscript's clarity. The authors reviewed the suggestions and take full responsibility for all mathematical results and proofs.

\end{document}